# Stabilization of Nonlinear Delay Systems Using Approximate Predictors and High-Gain Observers

Iasson Karafyllis[*a] and Miroslav Krstic[b]

[*] Corresponding author. Phone +30-210-7724478
[a] Dept. of Mathematics, National Technical University of Athens,15780, Athens, Greece, email: iasonkar@central.ntua.gr
[b] Dept. of Mechanical and Aerospace Eng., University of California, San Diego, La Jolla, CA 92093-0411, U.S.A., email: krstic@ucsd.edu

**Abstract**

We provide a solution to the heretofore open problem of stabilization of systems with arbitrarily long delays at the input and output of a nonlinear system using output feedback only. The solution is global, employs the predictor approach over the period that combines the input and output delays, addresses nonlinear systems with sampled measurements and with control applied using a zero-order hold, and requires that the sampling/holding periods be sufficiently short, though not necessarily constant. Our approach considers a class of globally Lipschitz strict-feedback systems with disturbances and employs an appropriately constructed successive approximation of the predictor map, a high-gain sampled-data observer, and a linear stabilizing feedback for the delay-free system. The obtained results guarantee robustness to perturbations of the sampling schedule and different sampling and holding periods are considered. The approach is specialized to linear systems, where the predictor is available explicitly.

## 1. Introduction

Summary of Results of the Paper. Even though numerous results have been developed in recent years for the stabilization of nonlinear systems with input delays by *state feedback* [15,17,19,20,21,22,25,26,37], and although additional delays in state measurements are allowed in our recent work [17], the problem of stabilization of systems with arbitrarily long delays at the input and/or output by *output feedback* has remained open.

In this work we provide a solution to this problem. Our solution addresses nonlinear systems with sampled measurements and with control applied using a zero-order hold, with a requirement that the sampling/holding periods be sufficiently short, though not necessarily constant. Our solution also employs the predictor approach to provide the control law with an estimate of the future state over a period that combines the input and output delays.

Our approach considers a class of globally Lipschitz strict-feedback systems with disturbances and employs an appropriately constructed successive approximation of the predictor map, a high-gain sampled-data observer, and a linear stabilizing feedback for the delay-free system. The obtained results can be applied to the linear time-invariant case as well, providing robust global sampled-data stabilizers, which are completely insensitive to perturbations of the sampling schedule and guarantee exponential convergence in the absence of measurement and modelling errors.

Our approach achieves input-to-state stability with respect to plant disturbances and measurement disturbances, as well as global exponential stability in the absence of disturbances.

Problem Statement and Literature. As in [15,17,19,20,21,22] we consider nonlinear systems of the form:
$$\dot{x}(t) = f(x(t), u(t-\tau), d(t))$$
$$x(t) \in \Re^n, u(t) \in \Re^m, d(t) \in \Re^l \tag{1.1}$$
where $\tau \geq 0$ is the input delay and $f: \Re^n \times \Re^m \times \Re^l \to \Re^n$ is a locally Lipschitz mapping with $f(0,0,0) = 0$. We employ the predictor-based approach, which is ubiquitous for linear systems (see the references in [20,21]) and is different from other approaches for systems with input delays [25,26,37], where the stabilizing feedback for the delay free system is either applied or is modified and stability is guaranteed for sufficiently small input delays. The input in (1.1) can be applied continuously or with zero-order hold (see [17]) and the measured output is usually assumed to be the state vector $x(t) \in \Re^n$. In [17], we extended predictor-based nonlinear control to the disturbance-free case (i.e, $d \equiv 0$) of sampled measurements and measurement delays expressed as
$$y(t) = x(\tau_i - r), \text{ for } t \in [\tau_i, \tau_{i+1})$$
where $y$ is the measured output, the discrete time instants $\tau_i$ are the sampling times and $r \geq 0$ is the measurement delay. The motivation is that sampling arises simultane-ously with input and output delays in control over networks. Few papers have studied this problem (exceptions are [13] where input and measurement delays are considered for linear systems but the measurement is not sampled and [18] where the unicycle is studied).

In the absence of delays, in sampled-data control of nonlinear systems semiglobal practical stability is generally guaranteed [8,29,30], with the desired region of attraction achieved by sufficiently fast sampling. Alternatively, global results are achieved under restrictive conditions on the structure of the system [7,12,32,39]. Simultaneous consideration to sampling and delays (either physical or sampling-induced) is given in the literature on control of linear and nonlinear systems over networks [5,6,11,30,32,35,36,39,41], but almost all available results rely on delay-dependent conditions for the existence of stabilizing feedback and in most cases the stability domain is depends on the sampling interval/ delay. Exceptions are the papers [2,23], where prediction-based control methodologies are employed.

The assumption that the state vector is measured is seldom realistic. Instead, measurement is a function of the state vector, i.e., the measured output of system (1.1) is given by:

$$y(t) = h(x(\tau_i - r)) + \xi(\tau_i), \quad t \in [\tau_i, \tau_{i+1}), i \in Z^+ \quad (1.2)$$

where $\{\tau_i\}_{i=0}^{\infty}$ is the set of sampling times being an increasing sequence with $\sup_{i \geq 0}(\tau_{i+1} - \tau_i) \leq T_1$, $T_1 > 0$ is the upper diameter of the sampling partition, $r \geq 0$ is the measurement delay, $h: \Re^n \to \Re^k$ is a continuous vector field with $h(0) = 0$ (the output map) and $\xi \in \Re^k$ is the measurement error. The measurements are obtained at discrete time instants (the sampling times).

We study the following problem in this paper: find a feedback law, which utilizes the sampled measurements and applies the input with zero-order hold, given by

$$u(t) = u_j, \quad t \in [jT_2, (j+1)T_2), j \in Z^+ \quad (1.3)$$

where $T_2 > 0$ is the holding period, such that the closed-loop system (1.1) with (1.2), (1.3) satisfies the Input-to-State stability (ISS) property from the inputs $(d, \xi) \in \Re^l \times \Re^k$ for all sampling partitions with $\sup_{i \geq 0}(\tau_{i+1} - \tau_i) \leq T_1$.

<u>Solution Provided in the Paper.</u> The above problem is considered for the case of constant delays $\tau \geq 0$, $r \geq 0$ and is solved for the class of globally Lipschitz systems of the form

$$\begin{aligned}\dot{x}_i(t) &= f_i(x_1(t),...,x_i(t)) + x_{i+1}(t) \\ &\quad + g_i(x(t),u(t))d_i(t), \quad i = 1,...,n-1 \\ \dot{x}_n(t) &= f_n(x(t)) + g_n(x(t),u(t))d_n(t) + u(t-\tau) \\ x(t) &= (x_1(t),...,x_n(t))' \in \Re^n, u(t) \in \Re, \\ d(t) &= (d_1(t),...,d_n(t))' \in \Re^n\end{aligned} \quad (1.4)$$

where $f_i : \Re^i \to \Re$ ($i = 1,...,n$) are globally Lipschitz functions with $f_i(0) = 0$ ($i = 1,...,n$) and the output map is $h(x) = x_1$. The inputs $d_i$ ($i = 1,...,n$) represent disturbances and the functions $g_i : \Re^i \to \Re$ ($i = 1,...,n$) are locally Lipschitz, bounded functions. In this case, we can show stabilizability of system (1.1) even under perturbations of the sampling schedule, by combining the sampled-data observer design in [14] and the approximate predictor control proposed in [15]. The feedback design is based on the corresponding delay free system

$$\begin{aligned}\dot{x}_i(t) &= f_i(x_1(t),...,x_i(t)) + x_{i+1}(t), \quad i = 1,...,n-1 \\ \dot{x}_n(t) &= f_n(x(t)) + u(t)\end{aligned} \quad (1.5)$$

The proposed control schemes for both cases consist of three components:

*1<sup>st</sup> Component:* An observer, which utilizes past input and output values in order to provide (continuous or discrete) estimates of the delayed state vector $x(t-r)$.

*2<sup>nd</sup> Component:* The predictor mapping that utilizes the estimation provided by the observer and past input values in order to provide an estimation of the future value of the state vector $x(t+\tau)$.

*3<sup>rd</sup> Component:* A nominal globally stabilizing feedback for the corresponding delay-free system.

The above control scheme has long been in use for linear systems [24,27,28,40] and it has been used even for partial differential equation systems [9], but is novel for nonlinear systems. Moreover, even for Linear Time-Invariant (LTI) systems

$$\dot{x}(t) = Ax(t) + Bu(t-\tau) + Gd(t) \quad (1.6)$$

where $x(t) \in \Re^n, u(t) \in \Re, d(t) \in \Re^n$, we provide new sampled-data feedback stabilizers that are robust to perturbations of the sampling schedule and guarantee exponential convergence in the absence of measurement and modeling errors.

<u>Notation.</u> We adopt the following notation:

* For a vector $x \in \Re^n$ we denote by $|x|$ its usual Euclidean norm, by $x'$ its transpose. For a real matrix $A \in \Re^{n \times m}$, $A' \in \Re^{m \times n}$ denotes its transpose and $|A| := \sup\{|Ax| ; x \in \Re^n, |x| = 1\}$ is its induced norm. $I \in \Re^{n \times n}$ denotes the identity matrix. By $A = \text{diag}(l_1, l_2,...,l_n)$ we mean a diagonal matrix with $l_1, l_2,...,l_n$ on its diagonal.

* $\Re_+$ denotes the set of non-negative real numbers. $Z^+$ denotes the set of non-negative integers. For every $t \geq 0$, $[t]$ denotes the integer part of $t \geq 0$, i.e., the largest integer being less or equal to $t \geq 0$. A partition $\pi = \{T_i\}_{i=0}^{\infty}$ of $\Re_+$ is an increasing sequence with $T_0 = 0$ and $T_i \to +\infty$.

* Let $x : [a-r, b) \to \Re^n$ with $b > a \geq 0$ and $r \geq 0$. By $T_r(t)x$ we denote the "history" of $x$ from $t-r$ to $t$, i.e., $(T_r(t)x)(\theta) := x(t+\theta); \theta \in [-r, 0]$, for $t \in [a,b)$. By

$\breve{T}_r(t)x$ we denote the "open history" of $x$ from $t-r$ to $t$, i.e., $(\breve{T}_r(t)x)(\theta) := x(t+\theta)$; $\theta \in [-r,0)$, for $t \in [a,b)$.

∗ Let $I \subseteq \Re$ be an interval. By $L^\infty(I;U)$ ($L^\infty_{loc}(I;U)$) we denote the space of measurable and (locally) bounded functions $u(\cdot)$ defined on $I$ and taking values in $U \subseteq \Re^m$. We do not identify functions in $L^\infty(I;U)$ which differ on a measure zero set. For $x \in L^\infty([-r,0];\Re^n)$ or $x \in L^\infty([-r,0);\Re^n)$ we define $\|x\|_r := \sup_{\theta \in [-r,0]} |x(\theta)|$ or $\|x\|_r := \sup_{\theta \in [-r,0)} |x(\theta)|$. The least upper bound $\sup_{\theta \in [-r,0]} |x(\theta)|$ is not the essential supremum but the actual supremum.

Throughout the paper, for $r=0$ we adopt the convention $L^\infty([-r,0];\Re^n) = \Re^n$ and $C^0([-r,0];\Re^n) = \Re^n$.

## 2. Globally Lipschitz Systems

We consider system (1.4) with output

$$y(\tau_i) = x_1(\tau_i - r) + \xi(\tau_i), \; i \in Z^+ \quad (2.1)$$

where $\{\tau_i\}_{i=0}^\infty$ is the set of sampling times and is a partition of $\Re_+$ with $\sup_{i \geq 0}(\tau_{i+1} - \tau_i) \leq T_1$. We assume that $r + \tau > 0$, where $r \geq 0$ is the measurement delay and $\tau \geq 0$ is the input delay. The locally bounded input $\xi : \Re_+ \to \Re$ represents the measurement error and the measurable and locally essentially bounded inputs $d_i : \Re_+ \to \Re$ ($i=1,...,n$) represent disturbances. Our main assumption is stated next.

**(A)** There exist constants $L \geq 0$ and $G \geq 0$ such that

$$|f_i(x_1,...,x_i) - f_i(z_1,...,z_i)| \leq L|(x_1-z_1,...,x_i-z_i)|,$$

$$\forall (x_1,...,x_i) \in \Re^i, \; \forall (z_1,...,z_i) \in \Re^i \quad (2.2)$$

$$|g_i(x,u)| \leq G, \; \forall(x,u) \in \Re^n \times \Re \quad (2.3)$$

for all $i=1,...,n$. Moreover, $f_i(0) = 0$ for all $i=1,...,n$.

Define $f(x) := (f_1(x_1),...f_n(x))' \in \Re^n$, $A = \{a_{i,j} : i,j=1,...n\} \in \Re^{n \times n}$ with $a_{i,i+1} = 1$ for all $i=1,...,n-1$ and $a_{i,j} = 0$ if $j \neq i+1$, $b = (0,...,0,1)' \in \Re^n$, $c := (1,0,...,0)' \in \Re^n$. Inequalities (2.2), (2.3) guarantee that system (1.4) is forward complete, i.e., for every $(x_0, u, d) \in \Re^n \times L^\infty_{loc}([-\tau, +\infty); \Re) \times L^\infty_{loc}(\Re_+; \Re^n)$ the solution $x(t) \in \Re^n$ of system (1.4) with initial condition $x(0) = x_0 \in \Re^n$ and corresponding to inputs $(u,d) \in L^\infty_{loc}([-\tau, +\infty); \Re) \times L^\infty_{loc}(\Re_+; \Re^n)$ exists for all $t \geq 0$. Indeed, the function $P(t) = |x(t)|^2/2$ satisfies

$$\dot{P}(t) \leq ((n+1)L+3)P(t) + G^2|d(t)|^2/2 + u^2(t-\tau)/2,$$

for almost all $t \geq 0$ for which the solution $x(t) \in \Re^n$ of system (1.4) exists. Integrating the previous differential inequality and using standard arguments, we conclude that the solution $x(t) \in \Re^n$ of system (1.4) exists for all $t \geq 0$ and satisfies the following estimate for all $t > 0$:

$$|x(t)| \leq \left(|x_0| + \frac{G \sup_{0 \leq s < t}|d(s)| + \sup_{-\tau \leq s < t-\tau}|u(s)|}{\sqrt{(n+1)L+3}}\right)\exp\left(\frac{(n+1)L+3}{2}t\right) \quad (2.4)$$

The proposed observer/predictor-based feedback law consists of three components:

1) A high-gain sampled-data observer for system (1.4), (2.1) which provides an estimate $z(t) \in \Re^n$ of the delayed state vector $x(t-r)$.

2) An approximate predictor, i.e., a mapping that utilizes the applied input values and the estimate $z(t) \in \Re^n$ provided by the observer in order to provide an estimate for $x(t+\tau)$.

3) A stabilizing feedback law for the delay-free system, i.e., system (1.5).

In what follows, we are going to describe the construction of each one of the components. We also assume that the input and measurement delay values $\tau, r \geq 0$ are perfectly known.

<u>1st Component (High-Gain Sampled-Data Observer):</u> Let $p = (p_1,...,p_n)' \in \Re^n$ be a vector such that the matrix $(A + pc') \in \Re^{n \times n}$ is Hurwitz. The existence of a vector $p = (p_1,...,p_n)' \in \Re^n$ is guaranteed, since the pair of matrices $(A,c)$ is observable. The proposed high-gain sampled-data observer is of the form:

$$\dot{z}_i(t) = f_i(z_1(t),...,z_i(t)) + z_{i+1}(t) + \theta^i p_i(c'z(t) - w(t)), i=1,...,n-1$$
$$\dot{z}_n(t) = f_n(z_1(t),...,z_n(t)) + \theta^n p_n(c'z(t) - w(t)) + u(t-r-\tau)$$
$$\dot{w}(t) = f_1(z_1(t)) + z_2(t) \quad, \quad t \in [\tau_i, \tau_{i+1})$$
$$w(\tau_{i+1}) = y(\tau_{i+1}) = x_1(\tau_{i+1} - r) + \xi(\tau_{i+1})$$
$$\tau_{i+1} = \tau_i + T_1 \exp(-b(\tau_i)), \tau_0 = 0$$

(2.5)

where $(z(t), w(t)) \in \Re^n \times \Re$, $\theta \geq 1$ is a constant to be chosen sufficiently large by the user and $b : \Re_+ \to \Re_+$ is an arbitrary non-negative locally bounded input that is unknown to the user. The sampling sequence $\{\tau_i\}_{i=0}^\infty$ is an arbitrary partition of $\Re_+$ with $\sup_{i \geq 0}(\tau_{i+1} - \tau_i) \leq T_1$, i.e., the sampling schedule is arbitrary. In order to justify the use of the high-gain sampled-data observer (2.5), we emphasize that system (2.5) is the feedback interconnection of the usual high-gain observer of system (1.4) which estimates

$x(t-r)$ instead of $x(t)$ and uses $w(t)$ instead of (the non-available signal) $x_1(t-r)$:

$$\dot{z}_i(t) = f_i(z_1(t),...,z_i(t)) + z_{i+1}(t) + \theta^i p_i(c'z(t)-w(t)), i=1,...,n-1$$
$$\dot{z}_n(t) = f_n(z_1(t),...,z_n(t)) + \theta^n p_n(c'z(t)-w(t)) + u(t-r-\tau)$$

and the inter-sample predictor of (the non-available signal) $x_1(t-r)$:

$$\dot{w}(t) = f_1(z_1(t)) + z_2(t) \quad , \quad t \in [\tau_i, \tau_{i+1})$$
$$w(\tau_{i+1}) = x_1(\tau_{i+1} - r) + \xi(\tau_{i+1})$$
$$\tau_{i+1} = \tau_i + T_1 \exp(-b(\tau_i)), \tau_0 = 0$$

which utilizes the measurements and predicts the value of $x_1(t-r)$ between two consecutive measurements. Sampled-data observers of this type (which are robust to sampling schedule perturbations) were proposed in [14,33,34].

<u>2$^{nd}$ Component (Approximate Predictor)</u>: Let $u \in L^\infty([0,T];\Re)$ be arbitrary and define the operator $P_{T,u} : C^0([0,T];\Re^n) \to C^0([0,T];\Re^n)$ by

$$(P_{T,u}x)(t) = x(0) + \int_0^t \left( f(x(\tau)) + Ax(\tau) + bu(\tau) \right) d\tau,$$
$$\text{for } t \in [0,T]. \quad (2.6)$$

We denote $P_{T,u}^l = \underbrace{P_{T,u} \ldots P_{T,u}}_{l \text{ times}}$ for every integer $l \geq 1$. We next define the operators $G_T : \Re^n \to C^0([0,T];\Re^n)$, $C_T : C^0([0,T];\Re^n) \to \Re^n$ and $Q_{T,u}^l : \Re^n \to \Re^n$ for $l \geq 1$ by

$$(G_T x_0)(t) = x_0, \text{ for } t \in [0,T] \text{ and } C_T x = x(T) \quad (2.7)$$
$$Q_{T,u}^l = C_T P_{T,u}^l G_T \quad (2.8)$$

We next define the mapping $P_{l,m}^u : \Re^n \to \Re^n$ for arbitrary $u \in L^\infty([0,r+\tau];\Re)$. Let $l,m \geq 1$ be integers and $T = (r+\tau)/m$. We define for all $x \in \Re^n$:

$$P_{l,m}^u x = Q_{T,u_m}^l \ldots Q_{T,u_1}^l x \quad (2.9)$$

where $u_i(s) = u(s+(i-1)T)$, $i=1,...,m$ for $s \in [0,T)$ ($u_i \in L^\infty([0,T];\Re)$ for $i=1,...,m$).

The operator $P_{l,m}^u$ is a nonlinear operator which provides an estimate of the value of the state vector of system (1.5) after $r+\tau$ time units when the input $u \in L^\infty([0,r+\tau];\Re)$ is applied. The operator is constructed based on the following procedure:
- first, we divide the time interval $[0,r+\tau]$ into $m \geq 1$ subintervals of equal length $T=(r+\tau)/m$,
- secondly, we apply the method of successive approximations to each one of the subintervals; more specifically we apply $l \geq 1$ successive approximations in order to get an estimate of the value of the state vector at the end of each one of the subintervals.

The following result was proved in [15] and is stated here for the convenience of the reader.

**Proposition 2.1 (see [15]):** *Let $l,m$ be positive integers with $(nL+1)T<1$, where $T=(r+\tau)/m$. Then there exists a constant $K := K(m) \geq 0$, independent of $l$, such that for every $u \in L^\infty([0,r+\tau];\Re)$ and $x \in \Re^n$ the following inequality holds:*

$$\left| P_{l,m}^u x - \phi(r+\tau,x;u) \right| \leq K \frac{((nL+1)T)^{l+1}}{1-(nL+1)T} \left( |x| + \sup_{0 \leq \tau < r+\tau} |u(\tau)| \right)$$
(2.10)

*where $\phi(t,x;u)$ denotes the unique solution of (1.5) at time $t \in [0,r+\tau]$, with initial condition $x \in \Re^n$ and corresponding to input $u \in L^\infty([0,r+\tau];\Re)$.*

Inequality (2.10) guarantees that by choosing $l,m$ sufficiently large then we can predict the value of the solution of (1.5) $r+\tau$ time units ahead, based only on the initial condition $x \in \Re^n$ and the applied input $u \in L^\infty([0,r+\tau];\Re)$. The prediction is given by $P_{l,m}^u x$.

Let $\delta_{r+\tau} : L_{loc}^\infty([-r-\tau,+\infty);\Re) \to L_{loc}^\infty([0,+\infty);\Re)$ denote the shift operator defined by

$$(\delta_{r+\tau}u)(t) := u(t-r-\tau), \text{ for } t \geq 0 \quad (2.11)$$

We are now able to define the approximate predictor mapping $\Phi_{l,m} : \Re^n \times L^\infty([-r-\tau,0];\Re) \to \Re^n$ defined by:

$$\Phi_{l,m}(x,u) := P_{l,m}^{\delta_{r+\tau}u} x \quad (2.12)$$

Using (2.2), (2.3), (2.10) and the Gronwall-Bellman lemma, we conclude that the following inequality holds for the solution of (1.4) for all $t \geq r$:

$$\left| \Phi_{l,m}(z,\breve{T}_{r+\tau}(t)u) - x(t+\tau) \right| \leq K \frac{((nL+1)T)^{l+1}}{1-(nL+1)T} \left( |z| + \sup_{t-r-\tau \leq s < t} |u(s)| \right)$$
$$+ \exp((nL+1)(r+\tau))(r+\tau)G \sup_{t-r \leq s \leq t+\tau} |d(s)|$$
$$+ \exp((nL+1)(r+\tau)) |z - x(t-r)|$$
(2.13)

More specifically, inequality (2.13) follows from (2.10) and the fact that

$$\left| \Phi_{l,m}(z,\breve{T}_{r+\tau}(t)u) - x(t+\tau) \right|$$
$$\leq \left| \Phi_{l,m}(z,\breve{T}_{r+\tau}(t)u) - \hat{x}(t+\tau) \right| + \left| \hat{x}(t+\tau) - x(t+\tau) \right|$$

where $\hat{x}(t)$ is the solution of (1.4) with initial condition $\hat{x}(t-r) = z$ corresponding to input $d \equiv 0$.

By virtue of (2.4) and (2.13), we obtain the following inequality for all $(u,z) \in L^\infty([-r-\tau,0];\Re) \times \Re^n$:

$$\left| \Phi_{l,m}(z,u) \right| \leq \Gamma \left( |z| + \sup_{-r-\tau \leq s < 0} |u(s)| \right) \quad (2.14)$$

where $\Gamma := K \frac{((nL+1)T)^{l+1}}{1-(nL+1)T} + \exp\left( \frac{(n+1)L+3}{2}(r+\tau) \right)$.

<u>3$^{rd}$ Component (Delay-Free Stabilizing Feedback)</u>: Due to the triangular structure of system (1.4), the results in [38] in conjunction with (2.2), (2.3), imply that there exists $k \in \Re^n$, a symmetric positive definite matrix $P \in \Re^{n \times n}$ and constants $\mu, \gamma > 0$ such that

$$x'P(A+bk')x + x'Pf(x) + x'P\,\mathrm{diag}(g_1(x,u),\ldots,g_n(x,u))d \leq -2\mu x'Px + \gamma |d|^2, \quad (2.15)$$

for all $(x,d,u) \in \Re^n \times \Re^n \times \Re$

We are now in a position to construct a stabilizing observer-based predictor feedback. Let $T_2 > 0$ be the "holding period". The feedback law is given by (2.5) with

$$u(t) = k'\Phi_{l,m}(z(iT_2), \breve{T}_{r+\tau}(iT_2)u),$$

for $t \in [iT_2, (i+1)T_2)$, $i \in Z^+$ (2.16)

where $\Phi_{l,m}(x,u)$ is defined by (2.12), (2.11), (2.9), (2.8), (2.7), (2.6) for integers $l, m \geq 1$.

In order to be able to show that the dynamic feedback law (2.5), (2.16) is successful, we need to assume that the upper diameter of the sampling partition and the holding period are sufficiently small. This is made in the following assumption.

**(B)** *Let* $Q \in \Re^{n \times n}$ *be a symmetric positive definite matrix that satisfies* $Q(A+pc') + (A'+cp')Q + 2qI \leq 0$ *for certain constant* $q > 0$ *and certain* $p \in \Re^n$. *Let* $P \in \Re^{n \times n}$ *be a symmetric positive definite matrix that satisfies (2.15) for certain constant* $\mu, \gamma > 0$ *and certain* $k \in \Re^n$. *The upper diameter of the sampling partition* $T_1 > 0$ *and the holding period* $T_2 > 0$ *are given or chosen by the user as sufficiently small so that the following inequalities hold:*

$$4|Qp|(L + \max(1, 2|Q|L\sqrt{n}/q))T_1\sqrt{|Q|/a} < q \quad (2.17)$$

$$\left((nL+1+|k|)\sqrt{\frac{b'Pb}{2K_1}} + \mu\right)|k|T_2 < \mu \quad (2.18)$$

*where* $a > 0$ *is a constant satisfying* $a|x|^2 \leq x'Qx$ *for all* $x \in \Re^n$, $0 < K_1$ *is a constant satisfying* $K_1|x|^2 \leq x'Px$ *for all* $x \in \Re^n$, $K > 0$ *is the constant involved in (2.13) and* $T = (r+\tau)/m$.

The following theorem guarantees that an appropriate selection of the parameters of the dynamic feedback law (2.5), (2.16) can guarantee the ISS property for the closed-loop system (1.4) with (2.5), (2.16).

**Theorem 2.2:** *Consider system (1.4) under assumptions (A), (B). Then for every* $\theta \geq 1$ *and for every pair of integers* $l, m > 0$ *chosen by the user so that* $m > (nL+1)(r+\tau)$ *and to satisfy the following inequalities*

$$4|Qp|(L+\theta)T_1\sqrt{|Q|/a} < q \quad (2.19)$$

$$\theta \geq \max(1, 2|Q|L\sqrt{n}/q) \quad (2.20)$$

$$\left((nL+1+|k|)\sqrt{\frac{b'Pb}{2K_1}} + \mu\right)|k|\left(T_2 + K\frac{((nL+1)T)^{l+1}}{1-(nL+1)T}\right) < \mu \quad (2.21)$$

*where* $a, K_1, K > 0$ *are the constants involved in assumption (B) and* $T = (r+\tau)/m$, *there exist constants* $\sigma > 0$, $\Theta_i > 0$ ($i = 1,\ldots,6$) *and a non-decreasing function* $M \in C^0(\Re_+; \Re_+)$, *such that for every* $x_0 \in C^0([-r,0]; \Re^n)$, $u_0 \in L^\infty([-r-\tau,0]; \Re)$, $(z_0, w_0) \in \Re^n \times \Re$, $(\xi, b, d) \in L^\infty_{loc}(\Re_+; \Re \times \Re_+ \times \Re^n)$ *the solution* $(T_r(t)x, \breve{T}_{r+\tau}(t)u, z(t), w(t))$ *of the closed-loop system (1.4), (2.5) and (2.16) with initial condition* $\breve{T}_{r+\tau}(0)u = u_0$, $T_r(0)x = x_0$, $(z(0), w(0)) = (z_0, w_0)$ *and corresponding to inputs* $(\xi, b, d) \in L^\infty_{loc}(\Re_+; \Re \times \Re_+ \times \Re^n)$ *satisfies the following inequality for all* $t \geq 0$:

$$|z(t)| + |w(t)| + \|T_r(t)x\|_r + \|\breve{T}_{r+\tau}(t)u\|_{r+\tau}$$
$$\leq \exp(-\sigma t) M\left(\sup_{0 \leq s \leq jT_2+\tau} b(s)\right)\left(\Theta_1|z_0| + \Theta_2|w_0| + \Theta_3\|x_0\|_r + \Theta_4\|u_0\|_{r+\tau}\right)$$
$$+ M\left(\sup_{0 \leq s \leq jT_2+\tau} b(s)\right)\Theta_5 \sup_{0 \leq s \leq t}\left(\exp(-\sigma(t-s))|\xi(s)|\right)$$
$$+ M\left(\sup_{0 \leq s \leq jT_2+\tau} b(s)\right)\Theta_6 \sup_{0 \leq s \leq t}\left(\exp(-\sigma(t-s))|d(s)|\right)$$
(2.22)

*where* $j = \min\{j \in Z^+ : jT_2 \geq r + T_1\}$.

By assumption (B), the user can select sufficiently large integers $l, m \geq 1$ so that inequality (2.21) holds. Indeed, the selection of sufficiently large integers $l, m \geq 1$ makes the term $C = K|k|\frac{((nL+1)T)^{l+1}}{1-(nL+1)T}$ sufficiently small: first we select an integer $m \geq 1$ so that $(nL+1)(r+\tau) < m$ and then (since $K := K(m) \geq 0$ is independent of $l \geq 1$; see Proposition 2.1) we can select a sufficiently large integer $l \geq 1$ so that $C$ becomes sufficiently small.

Clearly, inequality (2.22) is an ISS-like inequality, which guarantees the ISS property from the inputs $(\xi, d) \in L^\infty_{loc}(\Re_+; \Re \times \Re^n)$ in an almost uniform way for the input $b \in L^\infty_{loc}(\Re_+; \Re_+)$ for the closed-loop system (1.4), (2.5) and (2.16). More specifically, the effect of the inputs in (2.22) is expressed by means of "fading memory estimates" (see [16]), which are particularly useful for proving exponential stability in the case where $\xi$ or $d$ are functions of the state (for hybrid systems with delays the equivalence between "fading memory" estimates and "Sontag-like" estimates has not been established).

The proof of Theorem 2.2 is technical because the closed-loop system (1.4), (2.5) and (2.16) is a hybrid system which involves delays: for such systems even local existence of the solution is not trivial. The proof relies on the following methodology:

1) First, we prove that the solution of the closed-loop system (1.4), (2.5) and (2.16) exists for all times and for arbitrary initial conditions and inputs. This is achieved by Lemmas 2.3 and 2.4 below. Moreover, we show that the solution satisfies certain bounds which are useful for the subsequent analysis.

2) A second step (Lemma 2.5) is to show that the observer (2.5) provides estimates of the state vector which converge exponentially to the actual values of the state in the absence of errors.

3) A third step (Lemma 2.6) in the proof is to show that the applied control action (with Zero-Order Hold) is "close" to the control action that the nominal controller $u = k'x$ would give in the absence of input delays. However, in order to be able to guarantee this we have to require that sampling is fast enough and that the approximate predictor is sufficiently accurate.

4) Finally, the proof is completed by using all bounds that we have obtained in the previous steps and employing a small-gain argument.

The proofs of the following lemmas are provided in Appendix A.

**Lemma 2.3 (Bound on Observer State):** *Consider system (1.4) under the assumptions of Theorem 2.2. For every $x \in C^0([-r,+\infty);\Re^n)$, $u \in L^\infty_{loc}([-r-\tau,+\infty);\Re)$, $(z_0, w_0) \in \Re^n \times \Re$, $(\xi, b) \in L^\infty_{loc}(\Re_+; \Re \times \Re_+)$ the solution $(z(t), w(t)) \in \Re^n \times \Re$ of the hybrid system (2.5) with initial condition $(z(0), w(0)) = (z_0, w_0)$ and corresponding to inputs $(\xi, b) \in L^\infty_{loc}(\Re_+; \Re \times \Re_+)$, $(x, u)$ exists for all $t \geq 0$ and satisfies the following inequality:*

$$\exp(-2\omega t)\left(|z(t)|^2 + w^2(t)\right) \leq$$
$$|z_0|^2 + |w_0|^2 + \frac{1}{2\omega} \sup_{0 \leq s < t} |u(s-r-\tau)|^2 \qquad (2.23)$$
$$+ \frac{\left(\sup_{0 \leq s \leq t} |x(s-r)| + \sup_{0 \leq s \leq t} |\xi(s)|\right)^2}{1 - \exp\left(-2\omega T_1 \exp\left(-\sup_{0 \leq s \leq t} b(s)\right)\right)}$$

*where* $\omega := \max\left(L(n+1) + 2 + 2n \max_{i=1,\ldots,n}\left(\theta^{2i} p_i^2\right), 1 + L^2\right)/2$.

**Lemma 2.4 (Closed-Loop Solution Exists for all Times):** *Consider system (1.4) under the assumptions of Theorem 2.2. For every $x_0 \in C^0([-r,0];\Re^n)$, $u_0 \in L^\infty([-r-\tau,0);\Re)$, $(z_0, w_0) \in \Re^n \times \Re$, $(\xi, b, d) \in L^\infty_{loc}(\Re_+; \Re \times \Re_+ \times \Re^n)$ the solution $(T_r(t)x, \breve{T}_{r+\tau}(t)u, z(t), w(t))$ of the closed-loop system (1.4), (2.5) and (2.16) with initial condition $\breve{T}_{r+\tau}(0)u = u_0$, $T_r(0)x = x_0$, $(z(0), w(0)) = (z_0, w_0)$ and corresponding to inputs $(\xi, b, d) \in L^\infty_{loc}(\Re_+; \Re \times \Re_+ \times \Re^n)$ exists for all $t \geq 0$ and satisfies the following estimate:*

$$\sup_{0 \leq s \leq t}(|z(s)| + |w(s)|) + \sup_{-r \leq s \leq t}(|x(s)|) + \sup_{-r-\tau \leq s < t}(|u(s)|) \leq$$

$$\left(\frac{7(1+\Gamma)\exp(\beta T_2)}{\sqrt{1 - \exp\left(-2\omega T_1 \exp\left(-\sup_{0 \leq s \leq t} b(s)\right)\right)}}\right)^{g\left(\frac{t}{T_2}\right)} \Xi$$

$$\Xi := |z_0| + |w_0| + \|x_0\|_r + \|u_0\|_{r+\tau} + \sup_{0 \leq s \leq t}|\xi(s)| + G \sup_{0 \leq s \leq t}|d(s)|$$
(2.24)

*where* $g(t) := \min\{k \in Z^+ : t \leq k\}$, $\beta := \omega + \frac{(n+1)L+3}{2}$,

$\omega := \max\left(L(n+1) + 2 + 2n \max_{i=1,\ldots,n}\left(\theta^{2i} p_i^2\right), 1 + L^2\right)/2$ *and*

$\Gamma := K\left(\frac{((nL+1)T)^{l+1}}{1-(nL+1)T} + \exp\left(\frac{(n+1)L+3}{2}(r+\tau)\right)\right)$.

As remarked above, having completed the first step of the proof of Theorem 2.2 (which guarantees existence of the solution of the closed-loop system (1.4), (2.5) and (2.16) for all times and for arbitrary initial conditions and inputs), we are ready to proceed to the second step: to show that the observer (2.5) can provide estimates of the state vector. This is achieved by the following lemma.

**Lemma 2.5 (Convergence of Observer Estimate for Fast Sampling and High Observer Gain):** *Consider system (1.4) under the assumptions of Theorem 2.2. Then there exist constants $\sigma > 0$, $A_i > 0$ ($i = 1,\ldots,7$), which are independent of $T_2 > 0$ and $l, m$, such that for every $x_0 \in C^0([-r,0];\Re^n)$, $u_0 \in L^\infty([-r-\tau,0);\Re)$, $(z_0, w_0) \in \Re^n \times \Re$, $(\xi, b, d) \in L^\infty_{loc}(\Re_+; \Re \times \Re_+ \times \Re^n)$ the solution $(T_r(t)x, \breve{T}_{r+\tau}(t)u, z(t), w(t))$ of the closed-loop system (1.4), (2.5) and (2.16) with initial condition $\breve{T}_{r+\tau}(0)u = u_0$, $T_r(0)x = x_0$, $(z(0), w(0)) = (z_0, w_0)$ and corresponding to inputs $(\xi, b, d) \in L^\infty_{loc}(\Re_+; \Re \times \Re_+ \times \Re^n)$ satisfies the following estimate for all $t \geq r + T_1$:*

$$|z(t) - x(t-r)| \leq A_1 \exp(-\sigma(t-r))|z(r) - x(0)|$$
$$+ A_2 \sup_{0 \leq s \leq t}(\exp(-\sigma(t-s))|\xi(s)|)$$
$$+ A_3 \sup_{r \leq s \leq r+T_1}(\exp(-\sigma(t-s))|w(s) - x_1(s-r)|) \qquad (2.25)$$
$$+ A_4 \sup_{0 \leq s \leq t}(\exp(-\sigma(t-s))|d(s)|)$$

$$\sup_{r+T_1 \leq s \leq t}(\exp(\sigma s)|w(s) - x_1(s-r)|) \leq$$
$$A_5 \sup_{0 \leq s \leq t}(\exp(\sigma s)|\xi(s)|) + A_6|z(r) - x(0)|$$
$$+ \sup_{r \leq s \leq r+T_1}(\exp(\sigma s)|w(s) - x_1(s-r)|) \qquad (2.26)$$
$$+ A_7 \sup_{0 \leq s \leq t}(\exp(\sigma s)|d(s)|)$$

As explained above, the third step of the proof of Theorem 2.2 is to show that the applied control action (with Zero-Order Hold) is "close" to the control action that the nominal controller $u = k'x$ would give in the absence of input delays. This is achieved by the following lemma.

**Lemma 2.6 (Zero-Order Hold Control Close to Nominal Control if Sampling is Fast and Approximate Predictor is Accurate):** *Consider system (1.4) under the assumptions of Theorem 2.2. Define $j = \min\{j \in Z^+ : jT_2 \geq r + T_1\}$. Then for all sufficiently small $\sigma > 0$ and for all $(x_0, u_0, z_0, w_0) \in C^0([-r, 0]; \Re^n) \times L^\infty([-r-\tau, 0); \Re) \times \Re^n \times \Re$, $(\xi, b, d) \in L^\infty_{loc}(\Re_+; \Re \times \Re_+ \times \Re^n)$ (independent of $\sigma > 0$) the solution $(T_r(t)x, \tilde{T}_{r+\tau}(t)u, z(t), w(t))$ of the closed-loop system (1.4), (2.5) and (2.16) with initial condition $\tilde{T}_{r+\tau}(0)u = u_0$, $T_r(0)x = x_0$, $(z(0), w(0)) = (z_0, w_0)$ and corresponding to inputs $(\xi, b, d) \in L^\infty_{loc}(\Re_+; \Re \times \Re_+ \times \Re^n)$ satisfies the following estimate for all $t \geq jT_2 + \tau$:*

$$\eta \exp(\sigma t) |u(t-\tau) - k'x(t)|$$
$$\leq C|k|\Xi \exp(\sigma r) \sup_{jT_2 - r \leq s < jT_2 + \tau} \left( \exp(\sigma s) |u(s-\tau) - k'x(s)| \right)$$
$$+ A_1 \Xi \exp(\sigma r) |k|(C + \Omega) |z(r) - x(0)|$$
$$+ A_2 \Xi |k|(C + \Omega) \sup_{0 \leq s \leq t}(\exp(\sigma s)|\xi(s)|)$$
$$+ A_3 \Xi |k|(C + \Omega) \sup_{r \leq s \leq r + T_1} \left( \exp(\sigma s)|w(s) - x_1(s-r)| \right)$$
$$+ |k|\Xi [A_4 C + \Omega(A_4 + G(r+\tau)\exp(\sigma r)) + DG] \sup_{0 \leq s \leq t}(\exp(\sigma s)|d(s)|)$$
$$+ |k|\Xi [C(1+|k|)\exp(\sigma r) + D(nL+1+|k|)] \sup_{-r \leq s \leq t}(\exp(\sigma s)|x(s)|)$$

(2.27)

where $\Omega := \exp((nL+1)(r+\tau))$, $\Xi := \exp(\sigma(T_2 + \tau))$, $D := T_2 \exp(-\sigma\tau)$, $\eta := 1 - |k|T_2 - C|k|\exp(\sigma(T_2 + r + \tau))$, $C := K \dfrac{((nL+1)T)^{l+1}}{1 - (nL+1)T}$ and $A_i > 0$ ($i = 1, ..., 4$) are the constants involved in (2.25).

We are ready to give the proof of Theorem 2.2. The proof relies on the exploitation of inequalities (2.23), (2.24), (2.25) and (2.27) and use of a small-gain argument.

**Proof of Theorem 2.2:** Let $\sigma \in (0, \mu/2)$ be sufficiently small such that

$$\sqrt{\frac{b'Pb}{2K_1}} |k| \exp(\sigma T_2) [C(1+|k|)\exp(\sigma(r+\tau)) + T_2(nL+1+|k|)] < \eta \mu$$

and such that inequalities (2.25), (2.27) hold. The existence of sufficiently small $\sigma > 0$ satisfying $\sqrt{\dfrac{b'Pb}{2K_1}} |k| \exp(\sigma T_2) [C(1+|k|)\exp(\sigma(r+\tau)) + T_2(nL+1+|k|)] < \eta \mu$ is a direct consequence of (2.21). Define $V(t) = x'(t)Px(t)$. Using (2.15) we obtain the following differential inequality for almost all $t \geq 0$:

$$\dot{V}(t) \leq -2\mu V(t) + \frac{b'Pb}{2\mu}|u(t-\tau) - k'x(t)|^2 + 2\gamma|d(t)|^2 \quad (2.28)$$

Inequality (2.28) gives the following estimate for all $t > 0$:

$$|x(t)| \leq \sqrt{\frac{K_2}{K_1}} \exp(-\mu t)|x(0)|$$
$$+ \frac{1}{\mu}\sqrt{\frac{b'Pb}{2K_1}} \sup_{0 \leq s \leq t}\left(\exp\left(-\frac{\mu}{2}(t-s)\right)|u(s-\tau) - k'x(s)|\right) \quad (2.29)$$
$$+ \sqrt{\frac{2\gamma}{\mu K_1}} \sup_{0 \leq s \leq t}\left(\exp\left(-\frac{\mu}{2}(t-s)\right)|d(s)|\right)$$

where $0 < K_1 \leq K_2$ are constants satisfying $K_1 |x|^2 \leq x'Px \leq K_2 |x|^2$ for all $x \in \Re^n$. Using the inequality $\sigma \leq \mu/2$ we conclude that the following inequality holds for all $t > 0$:

$$|x(t)| \leq \sqrt{\frac{K_2}{K_1}} \exp(-\sigma t)|x(0)| + \sqrt{\frac{2\gamma}{\mu K_1}} \sup_{0 \leq s < t}(\exp(-\sigma(t-s))|d(s)|)$$
$$+ \frac{1}{\mu}\sqrt{\frac{b'Pb}{2K_1}} \sup_{0 \leq s < t}(\exp(-\sigma(t-s))|u(s-\tau) - k'x(s)|)$$

(2.30)

Inequality (2.30) implies the following inequality for all $t > 0$:

$$\sup_{0 \leq s \leq t}(\exp(\sigma s)|x(s)|) \leq \sqrt{\frac{K_2}{K_1}}|x(0)| + \sqrt{\frac{2\gamma}{\mu K_1}} \sup_{0 \leq s < t}(\exp(\sigma s)|d(s)|)$$
$$+ \frac{1}{\mu}\sqrt{\frac{b'Pb}{2K_1}} \sup_{0 \leq s < t}(\exp(\sigma s)|u(s-\tau) - k'x(s)|)$$

(2.31)

Combining (2.31) and (2.27), we obtain the following inequality for all $t \geq jT_2 + \tau$, where $j = \min\{j \in Z^+ : jT_2 \geq r + T_1\}$:

$$\mu \sup_{0 \leq s \leq t}(\exp(\sigma s)|x(s)|) \leq \mu\sqrt{\frac{K_2}{K_1}}|x(0)| + \sqrt{\frac{2\gamma\mu}{K_1}} \sup_{0 \leq s \leq t}(\exp(\sigma s)|d(s)|)$$
$$+ \bar{S}(1 + \eta^{-1}C|k|\exp(\sigma(T_2 + r + \tau))) \sup_{0 \leq s < jT_2 + \tau}(\exp(\sigma s)|u(s-\tau) - k'x(s)|)$$
$$+ \frac{\bar{S}}{\eta} A_1 \Xi \exp(\sigma r)|k|(C + \Omega)|z(r) - x(0)|$$
$$+ \frac{\bar{S}}{\eta} A_2 \Xi |k|(C + \Omega) \sup_{0 \leq s \leq t}(\exp(\sigma s)|\xi(s)|)$$
$$+ \frac{\bar{S}}{\eta} A_3 \Xi |k|(C + \Omega) \sup_{r \leq s \leq r + T_1}(\exp(\sigma s)|w(s) - x_1(s-r)|)$$
$$+ \frac{\bar{S}}{\eta}|k|\Xi[A_4 C + \Omega(A_4 + G(r+\tau)\exp(\sigma r)) + DG] \sup_{0 \leq s \leq t}(\exp(\sigma s)|d(s)|)$$
$$+ \frac{\bar{S}}{\eta}|k|\Xi[C(1+|k|)\exp(\sigma r) + D(nL+1+|k|)]\|x_0\|_r$$
$$+ \frac{\bar{S}}{\eta}|k|\Xi[C(1+|k|)\exp(\sigma r) + D(nL+1+|k|)] \sup_{0 \leq s \leq t}(\exp(\sigma s)|x(s)|)$$

where $\bar{S} := \sqrt{b'Pb/(2K_1)}$. It is clear from the above inequality that there exist constants $B_i > 0$ ($i = 1, ..., 6$) so that the following inequality holds for all $t \geq jT_2 + \tau$, where $j = \min\{j \in Z^+ : jT_2 \geq r + T_1\}$:

$$\sup_{0 \leq s \leq t}(\exp(\sigma s)|x(s)|) \leq B_1 \|x_0\|_r + B_2 \sup_{0 \leq s \leq t}(\exp(\sigma s)|d(s)|)$$
$$+ B_3 \sup_{0 \leq s < jT_2 + \tau}(\exp(\sigma s)|u(s - \tau) - k'x(s)|)$$
$$+ B_4 |z(r) - x(0)| + B_5 \sup_{0 \leq s \leq t}(\exp(\sigma s)|\xi(s)|) \quad (2.32)$$
$$+ B_6 \sup_{r \leq s \leq r + T_1}(\exp(\sigma s)|w(s) - x_1(s - r)|)$$

provided that

$$\overline{S}|k|\exp(\sigma T_2)[C(1+|k|)\exp(\sigma(r+\tau)) + T_2(nL+1+|k|)] < \eta\mu$$

where $\eta := 1 - |k|T_2 - C|k|\exp(\sigma(T_2 + r + \tau))$, $C := K\frac{((nL+1)T)^{l+1}}{1-(nL+1)T}$.

Combining inequalities (2.27), (2.32), (2.25), (2.24) and inequality (2.26), we obtain the existence of constants $\Theta_i > 0$ ($i = 1,...,6$) satisfying inequality (2.22). The proof is complete. ◁

## 3. Specialization to Linear Time Invariant Systems

For the LTI case (1.6), where the pair of matrices $A \in \Re^{n \times n}$, $B \in \Re^n$ is stabilizable and the output is

$$y(\tau_i) = c'x(\tau_i - r) + \xi(\tau_i), \quad i \in Z^+ \quad (3.1)$$

where $\{\tau_i\}_{i=0}^{\infty}$ is a partition of $\Re^+$ with $\sup_{i \geq 0}(\tau_{i+1} - \tau_i) \leq T_1$ and the pair of matrices $A \in \Re^{n \times n}$, $c \in \Re^n$ is a detectable pair, we apply the observer-based predictor stabilization scheme described in Section 2. There exist vectors $k \in \Re^n$ and $p \in \Re^n$ such that the matrices $A + Bk'$ and $A + pc'$ are Hurwitz matrices. Moreover, the predictor mapping is given by the expression

$$\Phi(x, u) := \exp(A(r + \tau))x + \int_{-r-\tau}^{0} \exp(-As)Bu(s)ds$$

The above prediction scheme is exact for the case $d \equiv 0$. Therefore, the following corollary can be proved in exactly the same way with Theorem 2.2.

**Corollary 3.1:** *Assume that there exist vectors $k \in \Re^n$, $p \in \Re^n$ such that the matrices $A + Bk'$, $A + pc'$ are Hurwitz matrices. For sufficiently small holding period $T_2 > 0$ and for sufficiently small sampling period $T_1 > 0$, there exist constants $\sigma > 0$, $\Theta_i > 0$ ($i = 1,...,7$) and a non-decreasing function $M \in C^0(\Re_+; \Re_+)$, such that for every $x_0 \in C^0([-r, 0]; \Re^n)$, $u_0 \in L^{\infty}([-r-\tau, 0); \Re)$, $(z_0, w_0) \in \Re^n \times \Re$, $(\xi, b, d) \in L^{\infty}_{loc}(\Re_+; \Re \times \Re_+ \times \Re^n)$ the solution $(T_r(t)x, \breve{T}_{r+\tau}(t)u, z(t), w(t))$ of the closed-loop system consisting of (1.6) with*

$$\dot{z}(t) = Az(t) + Bu(t - r - \tau) + p(c'z(t) - w(t))$$
$$\dot{w}(t) = c'Az(t) + c'Bu(t - r - \tau) \quad , \quad t \in [\tau_i, \tau_{i+1}) \quad (3.2)$$
$$w(\tau_{i+1}) = y(\tau_{i+1}) = c'x(\tau_{i+1} - r) + \xi(\tau_{i+1})$$
$$\tau_{i+1} = \tau_i + T_1 \exp(-b(\tau_i)), \tau_0 = 0$$

$$u(t) = k'\exp(A(r+\tau))z(iT_2)$$
$$+ \int_{-r-\tau}^{0} k'\exp(-As)Bu(iT_2 + s)ds \quad , \text{ for } t \in [iT_2, (i+1)T_2)$$

*and initial condition $\breve{T}_{r+\tau}(0)u = u_0$, $T_r(0)x = x_0$, $(z(0), w(0)) = (z_0, w_0) \in \Re^n \times \Re$ and corresponding to inputs $(\xi, b, d) \in L^{\infty}_{loc}(\Re_+; \Re \times \Re_+ \times \Re^n)$ satisfies inequality (2.22) for all $t \geq 0$, where $j = \min\{j \in Z^+ : jT_2 \geq r + T_1\}$.*

The advantage of the sampled-data feedback stabilizer (3.2) compared to other stabilizers for (1.6) (see for example [24]) is that the closed-loop system (1.6), (3.2) is completely insensitive to perturbations of the sampling schedule (this is guaranteed by inequality (2.22) and the fact that possible perturbations of the sampling schedule are quantified by means of the input $b \in L^{\infty}_{loc}(\Re_+; \Re_+)$).

## 4. Illustrative Example

In this section we consider the two dimensional system

$$\dot{x}_1(t) = f(x_1(t)) + x_2(t) + d(t), \quad \dot{x}_2(t) = u(t - \tau) \quad (4.1)$$

where $d(t) \in \Re$ and $f(x) = x^2 \text{sgn}(x)/\sqrt{1+x^2}$. For this function we have $\sup_{x \in \Re}|f'(x)| = 4\sqrt{2}/(3\sqrt{3})$ and consequently system (4.1) is of the form (1.4) and satisfies the global Lipschitz assumption made in Section 2. The one-dimensional, disturbance-free version of system (4.1) was studied in [15]. Here, we study system (4.1) with output available at discrete time instants:

$$y(t) = x_1(iT_1 - r), \text{ for } t \in [iT_1, (i+1)T_1), i \in Z^+ \quad (4.2)$$

where $T_1 > 0$ is the sampling period and $r \geq 0$ is the measurement delay. The input $u(t)$ is applied with zero-order hold with holding period $T_2 > 0$. Theorem 2.2 implies that there exist constants $\Theta_i > 0$ ($i = 1,...,5$) and $\sigma > 0$ such that for every $(x_0, u_0, z_0, w_0) \in C^0([-r, 0]; \Re^n) \times L^{\infty}([-r-\tau, 0); \Re) \times \Re^n \times \Re$ and $d \in L^{\infty}_{loc}(\Re_+; \Re)$ the solution $(T_r(t)x, \breve{T}_{r+\tau}(t)u, z(t), w(t))$ of the closed-loop system (4.1) with

$$\dot{z}_1(t) = f(z_1(t)) + z_2(t) - 3\theta(z_1(t) - w(t))$$
$$\dot{z}_2(t) = -3\theta^2(z_1(t) - w(t)) + u(t - r - \tau) \quad (4.3)$$
$$\dot{w}(t) = f(z_1(t)) + z_2(t) \quad , \quad t \in [iT_1, (i+1)T_1), i \in Z^+$$
$$w((i+1)T_1) = y((i+1)T_1) = x_1((i+1)T_1 - r)$$
$$u(t) = k'\Phi_{l,m}(z(iT_2), \breve{T}_{r+\tau}(iT_2)u), \text{ for } t \in [iT_2, (i+1)T_2) \quad (4.4)$$

where $l, m \geq 1$ are integers, the operator $\Phi_{l,m}: \Re^2 \times L^\infty([-r-\tau, 0); \Re) \to \Re^2$ is defined by (2.12), $k = -(15, 8)' \in \Re^2$ and initial condition $\breve{T}_{r+\tau}(0)u = u_0$, $T_r(0)x = x_0$, $(z(0), w(0)) = (z_0, w_0) \in \Re^n \times \Re$ satisfies the following inequality for all $t \geq 0$:

$$
\begin{aligned}
&|z(t)| + |w(t)| + \|T_r(t)x\|_r + \|\breve{T}_{r+\tau}(t)u\|_{r+\tau} \leq \\
&\exp(-\sigma t)\left(\Theta_1|z_0| + \Theta_2|w_0| + \Theta_3\|x_0\|_r + \Theta_4\|u_0\|_r\right) \quad (4.5) \\
&+ \Theta_5 \sup_{0 \leq s \leq t}\left(\exp(-\sigma(t-s))|d(s)|\right)
\end{aligned}
$$

provided that $l, m$ are sufficiently large positive integers, $\theta \geq 1$ is sufficiently large and the sampling period $T_1 > 0$ and holding period $T_2 > 0$ are sufficiently small. The closed-loop system (4.1), (4.3), (4.4) was tested numerically for $r = \tau = 1/4$. It was found that the selection

$$l=1, \ m=2, \ (X_1, X_2) = \left(z_1 + \frac{1}{4}(z_2 + f(z_1)), z_2 + \int_{-1/2}^{-1/4} u(s)ds\right)',$$

$$\Phi_{1,2}(z_1, z_2, u) = \left(X_1 + \frac{1}{4}(X_2 + f(X_1)), X_2 + \int_{-1/4}^{0} u(s)ds\right)'$$

$$\theta = 1, \ T_2 = 0.01, \ T_1 = 3T_2 = 0.03 \quad (4.6)$$

was appropriate in order to guarantee the ISS property from the input $d \in L^\infty_{loc}(\Re_+; \Re)$ for the closed-loop system. Figures 1, 2 and 3 show the time evolution of the state and the input for initial conditions $x_1(s) = x_2(s) = 1$ for $s \in [-1/4, 0]$, $u(s) = -2$ for $s \in [-1/2, 0)$ and $z_1(0) = z_2(0) = w(0) = 0$ for the disturbance-free case ($d(t) \equiv 0$) and for a sinusoidal disturbance ($d(t) = 0.5\sin(t)$). It is shown that all variables converge to zero for the disturbance-free case, while all variables ultimately follow an oscillation pattern for the case of external periodic forcing. The disturbance of amplitude 0.5 generates state oscillations whose amplitude is almost 2. This is the consequence of the limitation to the achievable disturbance attenuation performance that is caused by the presence of the significant dead time $r + \tau = 1/2$.

## 5. Concluding Remarks

We have expanded the applicability of delay-compensating stabilizing feedback to nonlinear systems where only output measurement is available and where such measurement is subject to long delays. Our designs employ either exact or approximate predictor maps. We perform state estimation using high-gain sampled-data observers. Our results guarantee ISS in the presence of disturbances for globally Lipschitz systems, provided the sampling/holding periods are sufficiently short.

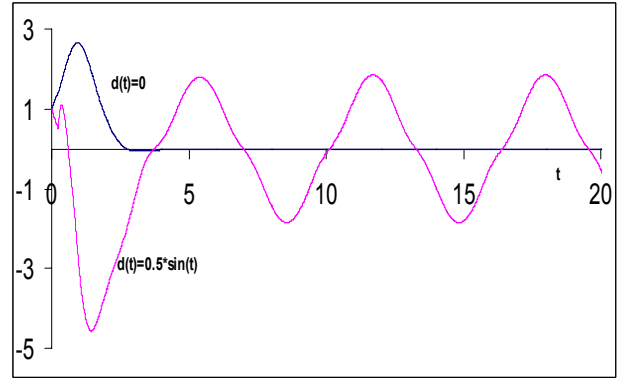

Fig. 1: Time evolution of the state $x_1(t)$ of the closed-loop system (4.1), (4.3), (4.4), (4.6)

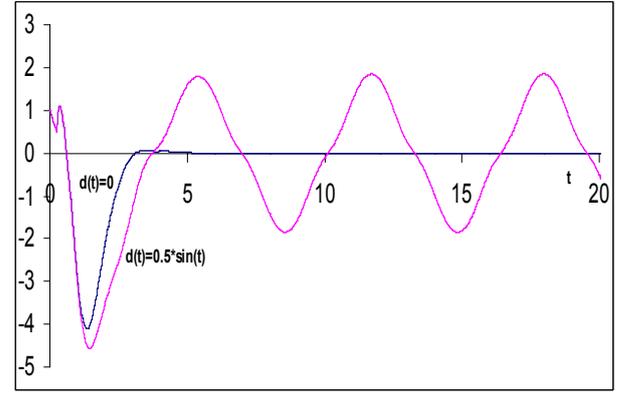

Fig. 2: Time evolution of the state $x_2(t)$ for the closed-loop system (4.1), (4.3), (4.4), (4.6)

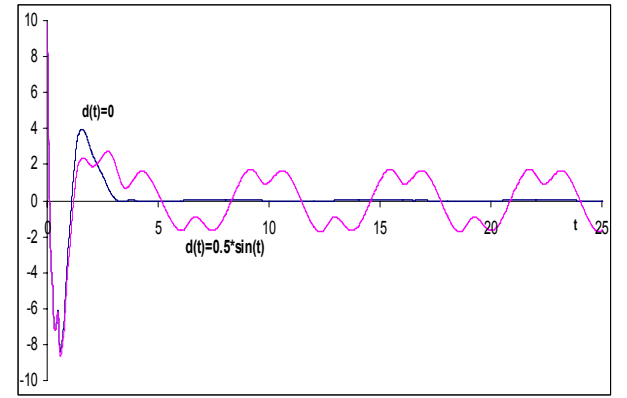

Fig. 3: Time evolution of the input $u(t)$ for the closed-loop system (4.1), (4.3), (4.4), (4.6)

Numerous relevant open problems remain that include multiple delays on inputs, states, and in the output map or quantization issues (as in [3,4,31]), or the possible use of emulation-based observers (as in [1]). Moreover, the issue of robustness with respect to variations of the input delay is crucial and can have serious effects (see for example [12]): it will be the topic of future work.

**References**

[1] Arcak, M. and D. Nešić, "A Framework for Nonlinear Sampled-Data Observer Design via Approximate Discrete-Time Models and Emulation", *Automatica*, 40, 2004, 1931-1938.

[2] Castillo-Toledo, B., S. Di Gennaro and G. Sandoval Castro, "Stability Analysis for a Class of Sampled Nonlinear Systems With Time-Delay", *Proceedings of the 49th Conference on Decision and Control*, Atlanta, GA, USA, 2010, 1575-1580.

[3] De Persis, C., "n-bit Stabilization of n-Dimensional Nonlinear Systems in Feedforward Form", *IEEE TAC*, 50(3), 2005, 285-297.

[4] De Persis, C., "Minimal Data-Rate Stabilization of Nonlinear Systems Over Networks With Large Delays", *Int. Journal of Robust and Nonlinear Control*, 20(10), 2010, 1097-1111.

[5] Fridman, E., A. Seuret and J.-P. Richard, "Robust Sampled-Data Stabilization of Linear Systems: An Input Delay Approach", *Automatica*, 40(8), 2004, 1441-1446.

[6] Gao, H., T. Chen and J. Lam, "A New System Approach to Network-Based Control", *Automatica*, 44(1), 2008, 39-52.

[7] Grüne, L., "Homogeneous State Feedback Stabilization of Homogenous Systems", *SIAM Journal on Control and Optimization*, 38(4), 2000, 1288-1308.

[8] Grüne, L. and D. Nešić, "Optimization Based Stabilization of Sampled-Data Nonlinear Systems via Their Approximate Discrete-Time Models", *SIAM Journal on Control and Optimization*, 42(1), 2003, 98-122.

[9] Guo, B.-Z., and K.-Y. Yang, "Output Feedback Stabilization of a One-Dimensional Schrödinger Equation by Boundary Observation With Time Delay", *IEEE TAC*, 55(5), 2010, 1226–1232.

[10] Heemels, M., A.R. Teel, N. van de Wouw and D. Nešić, "Networked Control Systems with Communication Constraints: Tradeoffs between Transmission Intervals, Delays and Performance", *IEEE TAC*, 55(8), 2010, 1781 - 1796.

[11] Herrmann, G., S.K. Spurgeon and C. Edwards, "Discretization of Sliding Mode Based Control Schemes", *Proceedings of the 38th Conference on Decision and Control*, Phoenix, Arizona, U.S.A., 1999, 4257-4262.

[12] Hetel, L., J. Daafouz, J. P. Richard and M. Jungers, "Delay-Dependent Sampled-Data Control Based on Delay Estimates", *Systems & Control Letters*, 60, 2011, 146–150.

[13] Jankovic, M., "Recursive Predictor Design for State and Output Feedback Controllers for Linear Time Delay Systems", *Automatica*, 46(3), 2010, 510-517.

[14] Karafyllis, I. and C. Kravaris, "From Continuous-Time Design to Sampled-Data Design of Observers", *IEEE Transactions on Automatic Control*, 54(9), 2009, 2169-2174.

[15] Karafyllis, I., "Stabilization By Means of Approximate Predictors for Systems with Delayed Input", *SIAM Journal on Control and Optimization*, 49(3), 2011, 1100-1123.

[16] Karafyllis, I., and Z.-P. Jiang, *Stability and Stabilization of Nonlinear Systems*, Springer-Verlag London (Series: Communications and Control Engineering), 2011.

[17] Karafyllis, I. and M. Krstic, "Nonlinear Stabilization under Sampled and Delayed Measurements, and with Inputs Subject to Delay and Zero-Order Hold", *IEEE Transactions on Automatic Control*, 57(5), 2012, 1141-1154.

[18] Kojima, K., T. Oguchi, A. Alvarez-Aguirre and H. Nijmeijer, "Predictor-Based Tracking Control of a Mobile Robot With Time-Delays", *Proceedings of NOLCOS 2010*, 167-172.

[19] Krstic, M., "Lyapunov tools for predictor feedbacks for delay systems: Inverse optimality and robustness to delay mismatch", *Automatica*, 44(11), 2008, 2930-2935.

[20] Krstic, M., *Delay Compensation for Nonlinear, Adaptive, and PDE Systems*, Birkhäuser Boston, 2009.

[21] Krstic, M., "Input Delay Compensation for Forward Complete and Strict-Feedforward Nonlinear Systems", *IEEE TAC*, 55(2), 2010, 287-303.

[22] Krstic, M., "Lyapunov Stability of Linear Predictor Feedback for Time-Varying Input Delay", *IEEE Transactions on Automatic Control*, 55(2), 2010, 554-559.

[23] Lozano, R., P. Castillo, P. Garcia and A. Dzul, "Robust Prediction-Based Control for Unstable Delay Systems: Application to the Yaw Control of a Mini-Helicopter", *Automatica*, 40(4), 2004, 603-612.

[24] Lozano, R., A. Sanchez, S. Salazar-Cruz and I. Fantoni, "Discrete-Time Stabilization of Integrators in Cascade: Real-Time Stabilization of a Mini-Rotorcraft" *International Journal of Control*, 81(6), 2008, 894–904.

[25] Mazenc, F., S. Mondie and R. Francisco, "Global Asymptotic Stabilization of Feedforward Systems with Delay at the Input", *IEEE TAC*, 49(5), 2004, 844-850.

[26] Mazenc, F., M. Malisoff and Z. Lin, "Further Results on Input-to-State Stability for Nonlinear Systems with Delayed Feedbacks", *Automatica*, 44(9), 2008, 2415-2421.

[27] Medvedev, A. and H. Toivonen, "Continuous-Time Deadbeat Observation Problem With Application to Predictive Control of Systems With Delay". *Kybernetika*, 30(6), 1994, 669–688.

[28] Mirkin, L., and N. Raskin, "Every Stabilizing Dead-Time Controller has an Observer–Predictor-Based Structure", *Automatica*, 39(10), 2003, 1747-1754.

[29] Nešić, D. and A.R. Teel, "Sampled-Data Control of Nonlinear Systems: An Overview of Recent Results", in Perspectives on Robust Control, R.S.O. Moheimani (Ed.), Springer-Verlag: New York, 2001, 221-239.

[30] Nešić, D. and A. Teel, "A Framework for Stabilization of Nonlinear Sampled-Data Systems Based on their Approximate Discrete-Time Models", *IEEE Transactions on Automatic Control*, 49(7), 2004, 1103-1122.

[31] Nešić, D. and D. Liberzon, "A unified framework for design and analysis of networked and quantized control systems", *IEEE TAC*, 54(4), 2009, 732-747.

[32] Nešić, D., A. R. Teel and D. Carnevale, "Explicit computation of the sampling period in emulation of controllers for nonlinear sampled-data systems", *IEEE TAC*, 54(3), 2009, 619-624.

[33] Postoyan, R., T. Ahmed-Ali and F. Lamnabhi-Lagarrigue, "Observers for classes of nonlinear networked systems", *Proceedings of the 6th IEEE. International Multi-Conference on Systems, Signals and Devices*, Djerba: Tunisia, 2009, 1-7.

[34] Postoyan, R., and D. Nesic, "A framework for the observer design for networked control systems", *Proceedings of the ACC*, Baltimore, U.S.A., 2010, 3678 - 3683.

[35] Tabbara, M., D. Nešić and A. R. Teel, "Networked control systems: emulation based design", in Networked Control Systems (Eds. D. Liu and F.-Y. Wang) Series in Intelligent Control and Intelligent Automation, World Scientific, 2007.

[36] Tabuada, P., "Event-Triggered Real-Time Scheduling of Stabilizing Control Tasks", *IEEE TAC*, 52(9), 2007, 1680-1685.

[37] Teel, A.R., "Connections between Razumikhin-Type Theorems and the ISS Nonlinear Small Gain Theorem", *IEEE TAC*, 43(7), 1998, 960-964.

[38] Tsinias, J., "A Theorem on Global Stabilization of Nonlinear Systems by Linear Feedback", *Systems and Control Letters*, 17(5), 1991, 357-362.

[39] Walsh, G. C., O. Beldiman, and L. G. Bushnell, "Asymptotic Behavior of Nonlinear Networked Control Systems", *IEEE Transactions on Automatic Control,*, 46(7), 2001, 1093-1097.

[40] Watanabe, K., and M. Sato, "A Predictor Control for Multivariable Systems With General Delays in Inputs and Outputs Subject to Unmeasurable Disturbances", *International Journal of Control*, 40(3), 1984, 435–448.

[41] Yu, M., L. Wang, T. Chu and F. Hao, "Stabilization of Networked Control Systems with Data Packet Dropout and Transmissions Delays: Continuous-Time Case", *European Journal of Control*, 11(1), 2005, 41-49.

**Appendix A. Supplementary material**

**Proof of Lemma 2.3:** Local existence and uniqueness follows from [16] (pages 23-27). Moreover, the analysis in [16] (pages 23-27) shows that the solution exists as long as it is bounded. In order to show that the solution remains bounded for all finite times, we consider the function $R(t) = |z(t)|^2/2 + w^2(t)/2$. By using algebraic manipulations and (2.2), it follows that the following differential inequality holds for almost all $t \in [\tau_i, \tau_{i+1})$ and $i \in Z^+$ for which the solution exists:

$$\dot{R}(t) \leq 2\omega R(t) + u^2(t-r-\tau)/2 \quad (A.1)$$

where $\omega := \max\left(L(n+1)+2+2n\max_{i=1,\ldots,n}\left(\theta^{2i}p_i^2\right), 1+L^2\right)/2$.

Integrating the differential inequality (A.1) we obtain for all $t \in [\tau_i, \tau_{i+1})$ and $i \in Z^+$ for which the solution exists:

$$|z(t)|^2 + w^2(t) \leq A\exp(2\omega(t-\tau_i))$$
$$A := |z(\tau_i)|^2 + |w(\tau_i)|^2 + \sup_{\tau_i \leq s < t}|u(s-r-\tau)|^2 \int_0^{t-\tau_i}\exp(-2\omega s)ds$$
(A.2)

Consequently, using a standard contradiction argument and (A.2), we are able to show that:

"if for some $i \in Z^+$ the solution exists at $t = \tau_i$ then the solution exists at $t = \tau_{i+1}$"

Using induction, (A.2) and the fact that $|w(\tau_i)| \leq \sup_{0 \leq s \leq \tau_i}|x(s-r)| + \sup_{0 \leq s \leq \tau_i}|\xi(s)|$ for all $i \in Z^+$ with $i \geq 1$, we show that the following inequality holds for all $i \in Z^+$ with $i \geq 2$:

$$\exp(-2\omega\tau_i)|z(\tau_i)|^2 \leq |z(0)|^2 + |w(0)|^2 + B$$
$$B = \left(\sup_{0 \leq s \leq \tau_i}|x(s-r)| + \sup_{0 \leq s \leq \tau_i}|\xi(s)|\right)^2 \sum_{k=1}^{i-1}\exp(-2\omega\tau_k)$$
$$+ \sup_{0 \leq s < \tau_i}|u(s-r-\tau)|^2 \int_0^{\tau_i}\exp(-2\omega s)ds$$
(A.3)

Inequalities (A.2), (A.3) and the fact that $|w(\tau_i)| \leq \sup_{0 \leq s \leq \tau_i}|x(s-r)| + \sup_{0 \leq s \leq \tau_i}|\xi(s)|$ for all $i \in Z^+$ with $i \geq 1$, show that the following inequality holds for all $i \in Z^+$ and $t \in [\tau_i, \tau_{i+1})$:

$$|z(t)|^2 + w^2(t) \leq \exp(2\omega t)\left(|z(0)|^2 + |w(0)|^2 + C\right)$$
$$C = \left(\sup_{0 \leq s \leq t}|x(s-r)| + \sup_{0 \leq s \leq t}|\xi(s)|\right)^2 \sum_{k=0}^{i}\exp(-2\omega\tau_k)$$
$$+ \sup_{0 \leq s < t}|u(s-r-\tau)|^2/(2\omega)$$
(A.4)

Inequality (2.23) is a direct consequence of (A.4) and the fact that $\tau_{i+1} \geq \tau_i + T_1\exp\left(-\sup_{0 \leq s \leq t}b(s)\right)$, which holds for all $i \in Z^+$ with $t \geq \tau_i$. The proof is complete. ◁

**Proof of Lemma 2.4:** We prove the lemma by proving the following claim for all $i \in Z^+$:

(Claim) For every $(z_0, w_0) \in \Re^n \times \Re$, $(x_0, u_0) \in C^0([-r,0];\Re^n) \times L^\infty([-r-\tau,0);\Re)$, $\xi \in L^\infty_{loc}(\Re_+;\Re)$, $(b,d) \in L^\infty_{loc}(\Re_+;\Re_+ \times \Re^n)$ the solution $(z(t), w(t)) \in \Re^n \times \Re$, $(T_r(t)x, \breve{T}_{r+\tau}(t)u) \in C^0([-r,0];\Re^n) \times L^\infty([-r-\tau,0);\Re)$ of the closed-loop system (1.4), (2.5) and (2.16) with initial condition $T_r(0)x = x_0$, $(z(0), w(0)) = (z_0, w_0) \in \Re^n \times \Re$, $\breve{T}_{r+\tau}(0)u = u_0$ and corresponding to inputs $(\xi, b, d) \in L^\infty_{loc}(\Re_+;\Re \times \Re_+ \times \Re^n)$ exists for all $t \in [0, iT_2]$ and satisfies (2.24) for all $t \in [0, iT_2]$.

It is clear that the claim holds for $i = 0$. Next assume that the claim holds for some $i \in Z^+$. Define

$$A_i := B \left( \frac{7(1+\Gamma)\exp(\beta T_2)}{\sqrt{1-\exp\left(-2\omega T_1 \exp\left(-\sup_{0\le s\le t} b(s)\right)\right)}} \right)^i$$

$$B = |z_0| + |w_0| + \|x_0\|_r + \|u_0\|_{r+\tau} + \sup_{0\le s\le t}|\xi(s)| + G\sup_{0\le s\le t}|d(s)|$$

Using (2.14), (2.16) and (2.24) for $t \in [0, iT_2]$, it is clear that $u(t)$ is well-defined on $[iT_2, (i+1)T_2)$ and satisfies the following inequality for all $t \in [iT_2, (i+1)T_2]$:

$$\sup_{-r-\tau \le s < t}(|u(s)|) \le \Gamma A_i \quad (A.5)$$

Using (2.4), (2.24) for $t \in [0, iT_2]$ and (A.5), it is clear that $x(t)$ is well-defined on $[iT_2, (i+1)T_2]$ and satisfies the following inequality for all $t \in [iT_2, (i+1)T_2]$:

$$|x(t)| \le \left((1+\Gamma)A_i + G\sup_{0\le s\le t}|d(s)|\right)\exp(((n+1)L+3)T_2/2) \quad (A.6)$$

Using Lemma 2.3, (A.5) and (A.6), it is clear that $(z(t), w(t))$ is well-defined on $[iT_2, (i+1)T_2]$ and satisfies the following inequality for all $t \in [iT_2, (i+1)T_2]$:

$$|z(t)| + |w(t)| \le 2C\exp(\omega T_2)\exp(((n+1)L+3)T_2/2)$$

$$C = \frac{G\sup_{0\le s\le t}|d(s)| + \sup_{0\le s\le t}|\xi(s)| + 2(1+\Gamma)A_i}{\sqrt{1-\exp\left(-2\omega T_1\exp\left(-\sup_{0\le s\le t}b(s)\right)\right)}} \quad (A.7)$$

Therefore, using the definition $\beta := \omega + ((n+1)L+3)/2$ and (A.5), (A.6), (A.7), we conclude that the following inequality holds for all $t \in [iT_2, (i+1)T_2]$:

$$\sup_{0\le s\le t}(|z(s)|+|w(s)|) + \sup_{-r\le s\le t}(|x(s)|) + \sup_{-r-\tau\le s<t}(|u(s)|)$$
$$\le \frac{7(1+\Gamma)\exp(\beta T_2)}{\sqrt{1-\exp\left(-2\omega T_1\exp\left(-\sup_{0\le s\le t}b(s)\right)\right)}} A_i$$

(A.8)

The fact that the claim holds for all $t \in [0, (i+1)T_2]$ is a direct consequence of (A.8). The proof is complete. ◁

**Proof of Lemma 2.5:** Define the quadratic error Lyapunov function $V(e) := e'\Delta_\theta^{-1}Q\Delta_\theta^{-1}e$, where $e(t) := z(t) - x(t-r)$, $\Delta_\theta := diag(\theta, \theta^2, ..., \theta^n)$. Using (2.2), (2.3), the identities $\Delta_\theta^{-1}A = \theta A\Delta_\theta^{-1}$, $c' = \theta c'\Delta_\theta^{-1}$ and the inequalities $\theta^{-i}|f_i(x_1+e_1,...,x_i+e_i) - f_i(x_1,...,x_i)| \le L|\Delta_\theta^{-1}e|$ for $i=1,...,n$ and all $(x,e) \in \Re^n \times \Re^n$ (which follow from (2.2)), we get for $\theta \ge \max(1, 2|Q|L\sqrt{n}/q)$ and for all $t \ge r$:

$$\dot{e}(t) = (A + \Delta_\theta pc')e(t) + \Delta_\theta p\eta(t) + \widetilde{p}(x(t-r), e(t))$$
$$- diag(g_1(x(t-r), u(t-r)),...,g_n(x(t-r), u(t-r)))d(t-r)$$
(A.9)

$$\dot{V}(t) \le -2\theta q|\Delta_\theta^{-1}e(t)|^2 + 2|\Delta_\theta^{-1}e(t)||Q||\Delta_\theta^{-1}\widetilde{p}(x(t-r), e(t))|$$
$$+ 2|e'(t)\Delta_\theta^{-1}||Qp||\eta(t)| + 2\theta^{-1}G|e'(t)\Delta_\theta^{-1}||Q||d(t-r)|$$
$$\le -2\theta q|\Delta_\theta^{-1}e(t)|^2 + 2|\Delta_\theta^{-1}e(t)|^2|Q|L\sqrt{n} + \theta q|e'(t)\Delta_\theta^{-1}|^2/2$$
$$+ 4\theta^{-1}q^{-1}|Qp|^2|\eta(t)|^2 + 4G^2\theta^{-3}q^{-1}|Q|^2|d(t-r)|^2$$
$$\le -\frac{\theta q}{2}|\Delta_\theta^{-1}e(t)|^2 + \frac{4}{\theta q}|Qp|^2|\eta(t)|^2 + \frac{4G^2}{\theta^3 q}|Q|^2|d(t-r)|^2$$
$$\le -\frac{\theta q}{2|Q|}V(e(t)) + \frac{4}{\theta q}|Qp|^2|\eta(t)|^2 + \frac{4G^2}{\theta^3 q}|Q|^2|d(t-r)|^2$$
(A.10)

where $\eta(t) = x_1(t-r) - w(t)$, $\widetilde{p}(x,e) = f(x+e) - f(x)$. Let $\sigma > 0$ be sufficiently small so that $4|Qp|(L+\theta)\exp(\sigma T_1)T_1\sqrt{|Q|/a} < q$ and $\sigma \le \theta q/(8|Q|)$. The existence of sufficiently small $\sigma > 0$ satisfying the inequality $4|Qp|(L+\theta)\exp(\sigma T_1)T_1\sqrt{|Q|/a} < q$ is guaranteed by (2.19). Using (A.10), we conclude that:

$$V(t) \le \exp(-4\sigma(t-r))V(r) +$$
$$16|Q|\theta^{-2}q^{-2}|Qp|^2\sup_{r\le s\le t}\left(\exp(-2\sigma(t-s))|\eta(s)|^2\right) +$$
$$16G^2\theta^{-4}q^{-2}|Q|^3\exp(2\sigma r)\sup_{0\le s\le t-r}\left(\exp(-2\sigma(t-s))|d(s)|^2\right)$$
(A.11)

for all $t \ge r$. Therefore, the following inequalities hold for all $t \ge r$:

$$|z(t) - x(t-r)| \le \exp(-2\sigma(t-r))\theta^{n-1}\sqrt{|Q|/a}|z(r) - x(0)| +$$
$$4|Qp|q^{-1}\theta^{n-1}\sqrt{|Q|/a}\sup_{r\le s\le t}\left(\exp(-\sigma(t-s))|\eta(s)|\right) +$$
$$4|Q|q^{-1}\theta^{n-2}G\sqrt{|Q|/a}\exp(\sigma r)\sup_{0\le s\le t-r}\left(\exp(-\sigma(t-s))|d(s)|\right)$$
(A.12)

$$|z_i(t) - x_i(t-r)| \leq \exp(-2\sigma(t-r))\theta^{i-1}\sqrt{|Q|/a}|z(r) - x(0)|$$
$$+ 4|Q|p|q^{-1}\theta^{i-1}\sqrt{|Q|/a}\sup_{r \leq s \leq t}(\exp(-\sigma(t-s))|\eta(s)|)$$
$$+ 4|Q|q^{-1}\theta^{i-2}G\sqrt{|Q|/a}\exp(\sigma r)\sup_{0 \leq s \leq t-r}(\exp(-\sigma(t-s))|d(s)|)$$
(A.13)

where $a > 0$ is a constant satisfying $a|x|^2 \leq x'Qx$ for all $x \in \Re^n$. Using (2.2) and (A.13), we obtain for $t \geq r$ a.e.:

$$|\dot{w}(t) - \dot{x}_1(t-r)| \leq \exp(-2\sigma(t-r))M|z(r) - x(0)| +$$
$$4|Q|p|q^{-1}M\sup_{r \leq s \leq t}(\exp(-\sigma(t-s))|\eta(s)|) + G|d_1(t-r)| +$$
$$4|Q|\theta^{-1}q^{-1}GM\exp(\sigma r)\sup_{0 \leq s \leq t-r}(\exp(-\sigma(t-s))|d(s)|)$$

where $M := (L+\theta)\sqrt{|Q|/a}$. The above inequality implies that the following estimate holds for all $t \in [\tau_i, \tau_{i+1})$, where $\tau_i$ with $i \geq 1$ is an arbitrary sampling time with $\tau_i \geq r$:

$$|\eta(t)| \leq |\xi(\tau_i)| + T_1 G \sup_{\tau_i - r \leq s \leq t-r}(|d_1(s)|)$$
$$+ \exp(-\sigma(\tau_i - r))T_1 M|z(r) - x(0)| +$$
$$4|Q|p|q^{-1}T_1 M \sup_{r \leq s \leq t}(\exp(-\sigma(\tau_i - s))|\eta(s)|)$$
$$+ 4|Q|\theta^{-1}q^{-1}GMT_1\exp(\sigma r)\sup_{0 \leq s \leq t-r}(\exp(-\sigma(\tau_i - s))|d(s)|)$$

Using the fact that $\tau_i \geq t - T_1$, the above inequalities give for all $t \in [\tau_i, \tau_{i+1})$, where $\tau_i$ with $i \geq 1$ is an arbitrary sampling time with $\tau_i \geq r$:

$$|\eta(t)| \leq \exp(\sigma T_1)\sup_{0 \leq s \leq t}(\exp(-\sigma(t-s))|\xi(s)|)$$
$$+ \exp(-\sigma(t-r))T_1 M \exp(\sigma T_1)|z(r) - x(0)|$$
$$+ 4|Q|p|q^{-1}\exp(\sigma T_1)T_1 M \sup_{r \leq s \leq t}(\exp(-\sigma(t-s))|\eta(s)|)$$
$$+ T_1 G \exp(\sigma(r+T_1))(4|Q|M\theta^{-1}q^{-1} + 1)\sup_{0 \leq s \leq t}(\exp(-\sigma(t-s))|d(s)|)$$
(A.14)

Notice that the above inequality holds for all $t \geq r + T_1$. It follows from (A.14) and the inequality $4|Q|p|(L+\theta)\exp(\sigma T_1)T_1\sqrt{|Q|/a} < q$ that the following inequality holds for all $t \geq r + T_1$:

$$\sup_{r+T_1 \leq s \leq t}(\exp(\sigma s)|\eta(s)|) \leq$$
$$\frac{q\exp(\sigma T_1)}{q - 4|Q|p|M\exp(\sigma T_1)T_1}\sup_{0 \leq s \leq t}(\exp(\sigma s)|\xi(s)|)$$
$$+ \frac{qM\exp(\sigma(r+T_1))T_1}{q - 4|Q|p|M\exp(\sigma T_1)T_1}|z(r) - x(0)|$$
$$+ \sup_{r \leq s \leq r+T_1}(\exp(\sigma s)|\eta(s)|)$$
$$+ \frac{(4|Q|M + q\theta)\exp(\sigma(r+T_1))}{\theta(q - 4|Q|p|M\exp(\sigma T_1)T_1)}GT_1\sup_{0 \leq s \leq t}(\exp(\sigma s)|d(s)|)$$
(A.15)

The existence of constants $\sigma > 0$, $A_i > 0$ ($i = 1,...,7$), which are independent of $T_2 > 0$ and $l, m$, satisfying (2.25) and (2.26) is a direct consequence of (A.12) and (A.15).

The proof is complete. ◁

**Proof of Lemma 2.6:** Let $\sigma > 0$ be sufficiently small such that (2.25) holds and such that $|k|T_2 + C|k|\exp(\sigma(T_2 + r + \tau)) < 1$, where $C := K\frac{((nL+1)T)^{l+1}}{1 - (nL+1)T}$. The existence of sufficiently small $\sigma > 0$ satisfying $|k|T_2 + C|k|\exp(\sigma(T_2 + r + \tau)) < 1$ is guaranteed by (2.21). Using (2.16), we obtain for all $i \in Z^+$ and $t \in [iT_2 + \tau, (i+1)T_2 + \tau)$:

$$|u(t - \tau) - k'x(t)| \leq |k||\Phi_{l,m}(z(iT_2), \breve{T}_{r+\tau}(iT_2)u) - x(t)|$$
$$\leq |k||\Phi_{l,m}(z(iT_2), \breve{T}_{r+\tau}(iT_2)u) - x(iT_2 + \tau)| + |k||x(iT_2 + \tau) - x(t)|$$
(A.16)

Using (2.13), we obtain for all $i \in Z^+$ with $iT_2 \geq r$ and $t \in [iT_2 + \tau, (i+1)T_2 + \tau)$:

$$|\Phi_{l,m}(z(iT_2), \breve{T}_{r+\tau}(iT_2)u) - x(iT_2 + \tau)|$$
$$\leq C|z(iT_2)| + C\sup_{iT_2 - r - \tau \leq s < iT_2}|u(s)|$$
$$+ \exp((nL+1)(r+\tau))G(r+\tau)\sup_{iT_2 - r \leq s \leq iT_2 + \tau}|d(s)|$$
$$+ \exp((nL+1)(r+\tau))|z(iT_2) - x(iT_2 - r)|$$
(A.17)

Combining (A.16) and (A.17) we obtain for all $i \in Z^+$ with $iT_2 \geq r$ and $t \in [iT_2 + \tau, (i+1)T_2 + \tau)$:

$$|u(t-\tau) - k'x(t)| \leq |k||x(iT_2 + \tau) - x(t)|$$
$$+ C|k| \sup_{iT_2 - r \leq s < iT_2 + \tau} |u(s-\tau) - k'x(s)|$$
$$+ |k|(C + \exp((nL+1)(r+\tau)))|z(iT_2) - x(iT_2 - r)| \quad (A.18)$$
$$+ G(r+\tau)|k|\exp((nL+1)(r+\tau)) \sup_{iT_2 - r \leq s \leq iT_2 + \tau} |d(s)|$$
$$+ C|k|(1+|k|) \sup_{iT_2 - r \leq s \leq iT_2 + \tau} |x(s)|$$

On the other hand, using (2.2) and (2.3), we conclude that the following inequality holds for all $i \in Z^+$ and $t \in [iT_2 + \tau, (i+1)T_2 + \tau)$:

$$\exp(\sigma t)|x(t) - x(iT_2 + \tau)|$$
$$\leq T_2(nL + 1 + |k|)\exp(\sigma T_2) \sup_{iT_2 + \tau \leq s \leq t}(\exp(\sigma s)|x(s)|) \quad (A.19)$$
$$+ T_2 \exp(\sigma t)|u(t-\tau) - k'x(t)|$$
$$+ T_2 G \exp(\sigma T_2) \sup_{0 \leq s \leq t}(\exp(\sigma s)|d(s)|)$$

Inequality (A.18) implies that the following inequality holds for all $i \in Z^+$ with $iT_2 \geq r$ and $t \in [iT_2 + \tau, (i+1)T_2 + \tau)$:

$$\exp(\sigma t)|u(t-\tau) - k'x(t)| \leq$$
$$C|k|\exp(\sigma(T_2 + r + \tau)) \sup_{iT_2 - r \leq s < iT_2 + \tau}(\exp(\sigma s)|u(s-\tau) - k'x(s)|)$$
$$+ \exp(\sigma t)|k|(C + \Omega)|z(iT_2) - x(iT_2 - r)|$$
$$+ \exp(\sigma t)|k||x(iT_2 + \tau) - x(t)|$$
$$+ G\Omega(r+\tau)|k|\exp(\sigma(T_2 + r + \tau)) \sup_{iT_2 - r \leq s < iT_2 + \tau}(\exp(\sigma s)|d(s)|)$$
$$+ C|k|(1+|k|)\exp(\sigma(T_2 + r + \tau)) \sup_{iT_2 - r \leq s < iT_2 + \tau}(\exp(\sigma s)|x(s)|)$$
$$(A.20)$$

where $\Omega := \exp((nL+1)(r+\tau))$. It follows from Lemma 2.5 and inequality (2.25) that the following inequality holds for all $i \in Z^+$ with $iT_2 \geq r + T_1$ and $t \in [iT_2 + \tau, (i+1)T_2 + \tau)$:

$$\exp(\sigma t)|z(iT_2) - x(iT_2 - r)| \leq A_1 \exp(\sigma(T_2 + \tau + r))|z(r) - x(0)|$$
$$+ A_2 \exp(\sigma(T_2 + \tau)) \sup_{0 \leq s \leq t}(\exp(\sigma s)|\xi(s)|)$$
$$+ A_3 \exp(\sigma(T_2 + \tau)) \sup_{r \leq s \leq r + T_1}(\exp(\sigma s)|w(s) - x_1(s-r)|)$$
$$+ A_4 \exp(\sigma(T_2 + \tau)) \sup_{0 \leq s \leq t}(\exp(\sigma s)|d(s)|)$$
$$(A.21)$$

Combining (A.20) and (A.21) we obtain for all $i \in Z^+$ with $iT_2 \geq r + T_1$ and $t \in [iT_2 + \tau, (i+1)T_2 + \tau)$:

$$\exp(\sigma t)(1 - |k|T_2)|u(t-\tau) - k'x(t)|$$
$$\leq C\Xi|k|\exp(\sigma r) \sup_{iT_2 - r \leq s < iT_2 + \tau}(\exp(\sigma s)|u(s-\tau) - k'x(s)|)$$
$$+ A_1\Xi\exp(\sigma r)|k|(C+\Omega)|z(r) - x(0)|$$
$$+ A_2\Xi|k|(C+\Omega) \sup_{0 \leq s \leq t}(\exp(\sigma s)|\xi(s)|)$$
$$+ A_3\Xi|k|(C+\Omega) \sup_{r \leq s \leq r + T_1}(\exp(\sigma s)|w(s) - x_1(s-r)|)$$
$$+ |k|\Xi[A_4 C + \Omega(A_4 + G(r+\tau)\exp(\sigma r)) + DG] \sup_{0 \leq s \leq t}(\exp(\sigma s)|d(s)|)$$
$$+ |k|\Xi[C(1+|k|)\exp(\sigma r) + D(nL+1+|k|)] \sup_{iT_2 - r \leq s \leq t}(\exp(\sigma s)|x(s)|)$$

where $\Xi := \exp(\sigma(T_2 + \tau))$ and $D := T_2 \exp(-\sigma\tau)$. Inequality (2.27) is a direct consequence of the above inequality. The proof is complete ◁

14